\documentclass[11pt]{article}
\usepackage{amsmath}
\usepackage{mathrsfs}
\usepackage{amsthm}
\usepackage{amsfonts}
\usepackage{amssymb}
\usepackage{bm}
\usepackage{mathtools}
\usepackage{color}

\usepackage{enumerate}
\usepackage{enumitem}
\setlist[enumerate,1]{label=(\arabic*),font=\textup,
leftmargin=7mm,labelsep=1.5mm,topsep=0mm,itemsep=-0.8mm}
\setlist[enumerate,2]{label=(\alph*).,font=\textup,
leftmargin=7mm,labelsep=1.5mm,topsep=-0.8mm,itemsep=-0.8mm}

\usepackage[pdfstartview=FitH,bookmarksopen=false,
colorlinks=true,linkcolor=blue,citecolor=blue]{hyperref}

\newtheorem{theorem}{Theorem}[section]
\newtheorem{definition}{Definition}[section]

\newtheorem{proposition}{Proposition}[section]
\newtheorem{question}{Question}[section]
\newtheorem{remark}{Remark}[section]

\title{\bf Sharp lower bounds on the spectral radius of uniform hypergraphs
concerning degrees\thanks {Research was partially supported by the National
Nature Science Foundation of China (grant numbers 11471210, 11571222)}}
\author {Lele Liu$^{1}$, \, Liying Kang$^{1}$, \, Erfang Shan$^{1,2}$
\thanks{\em Corresponding authors. Email address: efshan@shu.edu.cn (E. Shan),
lykang@i.shu.edu.cn (L. Kang), ahhylau@gmail.com (L. Liu)} \\
{\small $^{1}$Department of Mathematics, Shanghai University,
Shanghai 200444, P.R. China}\\
{\small$^{2}$School of Management, Shanghai University,
Shanghai 200444, P.R. China}}

\date{}
\begin{document}

\maketitle

\begin{abstract}

Let $\mathcal{A}(H)$ and $\mathcal{Q}(H)$ be the adjacency tensor and signless Laplacian
tensor of an $r$-uniform hypergraph $H$. Denote by $\rho(H)$ and $\rho(\mathcal{Q}(H))$
the spectral radii of $\mathcal{A}(H)$ and $\mathcal{Q}(H)$, respectively. In this paper,
we present a  lower bound on $\rho(H)$ in terms of vertex degrees and we characterize the
extremal hypergraphs attaining the bound, which solves a problem posed by Nikiforov
[V. Nikiforov, Analytic methods for uniform hypergraphs, Linear Algebra
Appl. 457 (2014) 455-535]. Also, we prove a lower bound on $\rho(\mathcal{Q}(H))$ concerning
degrees and give a characterization of the extremal
hypergraphs attaining the bound.

\bigskip \noindent{\bf Keywords:} Uniform hypergraph;
Adjacency tensor;
Signless Laplacian tensor;
Spectral radius

\medskip

\noindent{\bf AMS (2000) subject classification:}  15A42; 05C50; 05C70
\end{abstract}

\section{Introduction}
Let $G=(V(G),E(G))$ be a simple undirected graph with $n$ vertices, and $A(G)$
be the adjacency matrix of $G$. Let $\rho(G)$ be the spectral radius of $G$,
and $d_i$ be the degree of vertex $i$ of $G$, $i=1$, $2$, $\ldots$, $n$.
In 1988, Hofmeister \cite{Hofmeister} obtained a  lower bound on $\rho(G)$
in terms of degrees of vertices of $G$ as follows:
\begin{equation}\label{eq:Hofmeister}
\rho(G)\geqslant\left(\frac{1}{n}\sum_{i=1}^nd_i^2\right)^{\frac{1}{2}}.
\end{equation}
Furthermore, if $G$ is connected, then equality holds if and only if $G$ is either
a regular graph or a semiregular bipartite graph (see details in \cite{Hofmeister} and
\cite{AimeiYu}). The inequality \eqref{eq:Hofmeister} has many important applications
in spectral graph theory (see \cite{Nikiforov2010,Nikiforov2006,Nikiforov2007}).

In recent years the research on spectra of hypergraphs via tensors have
drawn increasingly extensive interest, accompanying with the rapid development
of tensor spectral theory. A {\em hypergraph} $H =(V,E)$  consists of a (finite)
set $V$ and a collection $E$ of non-empty
subsets of $V$ (see
\cite{Bretto}). The elements of $V$ are called {\em vertices} and the elements of
$E$ are called {\em hyperedges}, or simply {\em edges} of the
hypergraph.
If there is a risk of confusion we will denote the vertex set and the edge set of
a hypergraph $H$ explicitly
by $V(H)$ and $E(H)$, respectively.
An $r$-{\em uniform} hypergraph is a hypergraph
in which every edge has size $r$. Throughout this paper, we denote by $V(H)=[n]:=\{1,2, \ldots, n\}$
the vertex set of a hypergraph $H$.
For a vertex $i\in V(H)$, the {\em degree} of $i$, denoted
by $d_H(i)$ or simply by $d_i$, is the number
of edges containing $i$. If each vertex of $H$ has the same degree, we say that the
hypergraph $H$ is {\em regular}. For different $i$, $j\in V(H)$, $i$ and $j$ are said
to be {\em adjacent}, written $i\sim j$, if there is an edge of $H$ containing
both $i$ and $j$. A {\em walk} of hypergraph $H$ is defined to be an alternating
sequence of vertices and edges $i_1e_1i_2e_2\cdots i_{\ell}e_{\ell}i_{\ell+1}$
satisfying that $\{i_j,i_{j+1}\}\subseteq e_j\in E(H)$ for $1\leqslant j\leqslant \ell$.
A walk is called a {\em path} if all vertices and edges in the walk are distinct.
A hypergraph $H$ is called {\em connected} if for any vertices $i$, $j$, there
is a walk connecting $i$ and $j$. For positive integers $r$ and $n$, a real
{\em tensor} $\mathcal{A}=(a_{i_1i_2\cdots i_r})$ of order $r$ and dimension $n$
refers to a multidimensional array (also called {\em hypermatrix}) with entries
$a_{i_1i_2\cdots i_r}$ such that $a_{i_1i_2\cdots i_r}\in\mathbb{R}$ for
all $i_1$, $i_2$, $\ldots$, $i_r\in[n]$. We say that tensor $\mathcal{A}$
is {\em symmetric} if its entries $a_{i_1i_2\cdots i_r}$ are invariant
under any permutation of its indices.

Recently, Nikiforov \cite{Nikiforov} presented some analytic methods for studying
uniform hypergraphs, and posed the following question (see \cite[Question 11.5]{Nikiforov}):
\begin{question}[\cite{Nikiforov}]
\label{question}
Suppose that $H$ is an $r$-uniform hypergraph on $n$ vertices $(r\geqslant 3)$. Let $d_i$
be the degree of vertex $i$, $i\in[n]$, and $\rho(H)$ be the spectral radius of $H$. Is
it always true
\[
\rho(H)\geqslant\left(\frac1n\sum_{i=1}^nd^{\frac{r}{r-1}}_i\right)^{\frac{r-1}{r}}?
\]
\end{question}

In this paper, we focus on the above question, and give a solution to Question
\ref{question}. Our main results can be stated as follows.
\begin{theorem}
\label{thm:Main result}
Suppose that $H$ is an $r$-uniform hypergraph on $n$ vertices $(r\geqslant 3)$. Let $d_i$
be the degree of vertex $i$ of $H$, and $\rho(H)$ be the spectral radius of $H$. Then
\[
\rho(H)\geqslant
\left(\frac1n\sum_{i=1}^nd^{\frac{r}{r-1}}_i\right)^{\frac{r-1}{r}}.
\]
If $H$ is connected, then the equality holds if and only if $H$ is regular.
\end{theorem}

\begin{theorem}
\label{Signless Laplacian}
Let $H$ be a connected $r$-uniform hypergraph on $n$ vertices $(r\geqslant 3)$. Suppose
that $d_i$, $i\in[n]$, is the degree of vertex $i$, and $\rho(\mathcal{Q}(H))$ is the
spectral radius of the signless Laplacian tensor $\mathcal{Q}(H)$. Then
\[
\rho(\mathcal{Q}(H))\geqslant 2\left(\frac1n\sum_{i=1}^nd^{\frac{r}{r-1}}_i\right)^{\frac{r-1}{r}},
\]
with equality  if and only if $H$ is regular.
\end{theorem}
\section{Preliminaries}
\label{sec2}
In this section we review some basic notations and necessary conclusions.
Denote the set of nonnegative vectors (positive vectors) of dimension $n$
by $\mathbb{R}_+^n$ ($\mathbb{R}_{++}^n$). The unit tensor of order $r$
and dimension $n$ is the tensor $\mathcal{I}_n=(\delta_{i_1i_2\cdots i_r})$,
whose entry is $1$ if $i_1=i_2=\cdots=i_r$ and $0$ otherwise.

The following general product of tensors was defined by Shao
\cite{Shao:General product}, which is a generalization of the matrix case.

\begin{definition}[\cite{Shao:General product}]
Let $\mathcal{A}$ (and $\mathcal{B}$) be an order $r\geqslant 2$
(and order $k\geqslant 1$), dimension $n$ tensor. Define the product
$\mathcal{AB}$ to be the following tensor $\mathcal{C}$ of order
$(r-1)(k-1)+1$ and dimension $n$
\[
c_{i\alpha_1\cdots\alpha_{r-1}}=\sum_{i_2,\ldots,i_r=1}^na_{ii_2\cdots i_r}
b_{i_2\alpha_1}\cdots b_{i_r\alpha_{r-1}}
~~(i\in [n], \alpha_1,\ldots,\alpha_{r-1}\in [n]^{k-1}).
\]
\end{definition}

From the above definition, let $\bm{x}=(x_1,x_2,\ldots,x_n)^{\mathrm{T}}$ be a column
vector of dimension $n$. Then $\mathcal{A}\bm{x}$ is a vector in $\mathbb{C}^n$, whose
$i$-{\em th} component is as the following
\[
(\mathcal{A}\bm{x})_i=\sum_{i_2,\ldots,i_r=1}^na_{ii_2\cdots i_r}x_{i_2}\cdots x_{i_r},~~
i\in [n]
\]
and
\[
\bm{x}^{\mathrm{T}}(\mathcal{A}\bm{x})=\sum_{i_1,i_2,\ldots,i_r=1}^n
a_{i_1i_2\cdots i_r}x_{i_1}x_{i_2}\cdots x_{i_r}.
\]

In 2005, Lim \cite{Lim} and Qi \cite{Qi2005} independently introduced the concepts of
tensor eigenvalues and the spectra of tensors. Let $\mathcal{A}$ be an order $r$ and
dimension $n$ tensor, $\bm{x}=(x_1,x_2\ldots,x_n)^{\mathrm{T}}\in\mathbb{C}^n$
be a column vector of dimension $n$. If there exists a number $\lambda\in\mathbb{C}$
and a nonzero vector $\bm{x}\in\mathbb{C}^{n}$ such that
\[
\mathcal{A}\bm{x}=\lambda \bm{x}^{[r-1]},
\]
then $\lambda$ is called an {\em eigenvalue} of $\mathcal{A}$, $\bm{x}$ is called
an {\em eigenvector} of $\mathcal{A}$ corresponding to the eigenvalue $\lambda$,
where $\bm{x}^{[r-1]}=(x_1^{r-1},x_2^{r-1},\ldots,x_n^{r-1})^{\mathrm{T}}$. The
{\em spectral radius} $\rho(\mathcal{A})$ of $\mathcal{A}$ is the maximum modulus
of the eigenvalues of $\mathcal{A}$. It was proved that $\lambda$ is an eigenvalue
of $\mathcal{A}$ if and only if it is a root of the characteristic polynomial of
$\mathcal{A}$ (see details in \cite{Shao:Connected odd-bipartite}).

In 2012, Cooper and Dutle \cite{Cooper:Spectra Uniform Hypergraphs} defined the
adjacency tensors for $r$-uniform hypergraphs.
\begin{definition}
[\cite{Cooper:Spectra Uniform Hypergraphs,Qi2014}]
Let $H=(V(H),E(H))$ be an $r$-uniform hypergraph on $n$ vertices. The adjacency
tensor of $H$ is defined as the order $r$ and dimension $n$ tensor
$\mathcal{A}(H)=(a_{i_1i_2\cdots i_r})$, whose $(i_1i_2\cdots i_r)$-entry is
\[
a_{i_1i_2\cdots i_r}=\begin{cases}
\frac{1}{(r-1)!}, & \text{if}~\{i_1,i_2,\ldots,i_r\}\in E(H),\\
0, & \text{otherwise}.
\end{cases}
\]
Let $\mathcal{D}(H)$ be an order $r$ and dimension $n$ diagonal tensor with its
diagonal element $d_{ii\cdots i}$ being $d_i$, the degree of vertex $i$, for all
$i\in[n]$. Then $\mathcal{L}(H)=\mathcal{D}(H)-\mathcal{A}(H)$ is the Laplacian
tensor of $H$, and $\mathcal{Q}(H)=\mathcal{D}(H)+\mathcal{A}(H)$ is the signless
Laplacian tensor of $H$.
\end{definition}
For an $r$-uniform hypergraph $H$, denote the spectral radius of $\mathcal{A}(H)$
by $\rho(H)$. It should be announced that
spectral radius defined in \cite{Nikiforov} differ from this paper, while for an
$r$-uniform hypergraph $H$ the spectral radius defined in \cite{Nikiforov} equals
to $(r-1)!\rho(H)$. This is not essential and does not effect the result.

In \cite{Friedland}, the weak irreducibility of nonnegative tensors was defined. It
was proved that an $r$-uniform hypergraph $H$ is connected if and only if its adjacency
tensor $\mathcal{A}(H)$ is weakly irreducible (see \cite{Friedland} and
\cite{Yang:Nonegative Weakly Irreducible Tensors}). Clearly, this shows that if $H$
is connected, then $\mathcal{A}(H)$, $\mathcal{L}(H)$ and $\mathcal{Q}(H)$ are all
weakly irreducible.
The following result for nonnegative tensors is stated as a part of Perron-Frobenius theorem in \cite{K.C.Chang.etc:Perron-Frobenius Theorem}.
\begin{theorem}[\cite{K.C.Chang.etc:Perron-Frobenius Theorem}]
\label{thm:Perron-Frobenius}
Let $\mathcal{A}$ be a nonnegative tensor of order $r$ and dimension $n$. Then we have
the following statements.
\begin{enumerate}
\item $\rho(\mathcal{A})$ is an eigenvalue of $\mathcal{A}$ with a nonnegative
eigenvector corresponding to it.

\item If $\mathcal{A}$ is weakly irreducible, then $\rho(\mathcal{A})$ is the
unique eigenvalue of $\mathcal{A}$ with the unique eigenvector $\bm{x}\in\mathbb{R}_{++}^n$,
up to a positive scaling coefficient.
\end{enumerate}
\end{theorem}

\begin{theorem}
[\cite{Qi2013}]
\label{relaigh}
Let $\mathcal{A}$ be a nonnegative symmetric tensor of order $r$ and dimension $n$.
Then we have
\[
\rho(\mathcal{A})=\max\left\{\bm{x}^{\mathrm{T}}(\mathcal{A}\bm{x})\,|\, x\in\mathbb{R}_{+}^{n},
||\bm{x}||_r=1\right\}.
\]
Furthermore, $\bm{x}\in\mathbb{R}_{+}^{n}$ with $||\bm{x}||_r=1$ is an optimal
solution of the above optimization problem if and only if it is an eigenvector of
$\mathcal{A}$ corresponding to the eigenvalue $\rho(\mathcal{A})$.
\end{theorem}

The following concept of {\em direct products} (also called {\em Kronecker product}) of tensors
was defined in \cite{Shao:General product}, which is a generalization of the direct products
of matrices.
\begin{definition}[\cite{Shao:General product}]\label{defn:Direct product}
Let $\mathcal{A}$ and $\mathcal{B}$ be two order $r$ tensors with dimension $n$ and $m$,
respectively. Define the direct product $\mathcal{A}\otimes\mathcal{B}$ to be the following
tensor of order $r$ and dimension $mn$ (the set of subscripts is taken as $[n]\times [m]$
in the lexicographic order):
\[
(\mathcal{A}\otimes\mathcal{B})_{(i_1,j_1)(i_2,j_2)\cdots(i_r,j_r)}
=a_{i_1i_2\cdots i_r}b_{j_1j_2\cdots j_r}.
\]
In particular, if $\bm{x}=(x_1,x_2,\ldots,x_n)^{\mathrm{T}}$ and
$\bm{y}=(y_1,y_2,\ldots,y_m)^{\mathrm{T}}$ are two column vectors
with dimension $n$ and $m$, respectively. Then
\[
\bm{x}\otimes\bm{y}=(x_1y_1,\ldots,x_1y_m,
x_2y_1,\ldots,x_2y_m,\ldots,
x_ny_1,\ldots,x_ny_m)^{\mathrm{T}}.
\]
\end{definition}

The following basic results can be found in \cite{Shao:General product}.
\begin{proposition}[\cite{Shao:General product}]
\label{prop}
The following conclusions hold.
\begin{enumerate}
\item $(\mathcal{A}_1+\mathcal{A}_2)\otimes\mathcal{B}=
\mathcal{A}_1\otimes\mathcal{B}+\mathcal{A}_2\otimes\mathcal{B}$.
\item $(\lambda\mathcal{A})\otimes\mathcal{B}=\mathcal{A}\otimes(\lambda\mathcal{B})=
\lambda(\mathcal{A}\otimes\mathcal{B})$, $\lambda\in\mathbb{C}$.
\item $(\mathcal{A}\otimes\mathcal{B})(\mathcal{C}\otimes\mathcal{D})=(\mathcal{AC})\otimes(\mathcal{BD})$.
\end{enumerate}
\end{proposition}

\begin{theorem}[\cite{Shao:General product}]
\label{thm:Spectral H times G}
Suppose that $\mathcal{A}$ and $\mathcal{B}$ are two order $r$ tensors with dimension
$n$ and $m$, respectively. Let $\lambda$ be an eigenvalue of $\mathcal{A}$ with corresponding
eigenvector $\bm{u}$, and $\mu$ be an eigenvalue of $\mathcal{B}$ with corresponding
eigenvector $\bm{v}$. Then $\lambda\mu$ is an eigenvalue of $\mathcal{A}\otimes\mathcal{B}$
with corresponding eigenvector $\bm{u}\otimes \bm{v}$.
\end{theorem}

\section{Proof of Theorem 1.1}
In this section, we shall give a proof of Theorem \ref{thm:Main result}.

\noindent{\bfseries Proof of Theorem \ref{thm:Main result}}.
Let $H$ be an $r$-uniform hypergraph with spectral radius $\rho(H)$ and vertex
set $V(H)=[n]$, and denote by $d_i$ the degree of vertex $i$ of $H$, $i=1$, $2$, $\ldots$, $n$.

We now define an $r$-uniform hypergraph $\widetilde{H}$ as follows. Hypergraph $\widetilde{H}$
has vertex set $V(H)\times [r]$, and $\{(i_1,j_1),(i_2,j_2),\ldots,(i_r,j_r)\}\in E(\widetilde{H})$
is an edge of $\widetilde{H}$ if and only if $\{i_1,i_2,\ldots,i_r\}\in E(H)$ and $j_1,j_2,\ldots,j_r$ are
distinct each other. Let $\mathcal{A}(H)=(a_{i_1i_2\cdots i_r})$ be the adjacency tensor of $H$. We define
an order $r$ and dimension $r$ tensor $\mathcal{B}=(b_{j_1j_2\cdots j_r})$ as follows:
\begin{equation}
\label{eq:B}
b_{j_1j_2\cdots j_r}=
\begin{cases}
1,  & \text{if} ~j_1,j_2,\ldots,j_r~\text{are distinct each other},\\
0, & \text{otherwise}.
\end{cases}
\end{equation}
We denote the adjacency tensor of $\widetilde{H}$ by $\mathcal{A}(\widetilde{H})$,
in which the set of subscripts is taken as $[n]\times [r]$ in the lexicographic
order.

\noindent{\bfseries Claim 1.} If $H$ is connected, then $\widetilde{H}$ is connected.

\noindent{\bfseries Proof of Claim 1.}\,
It suffices to show that for any $(i,j)$, $(s,t)\in V(\widetilde{H})$, there exists a
walk connecting them. We distinguish the following two cases.

\noindent{\bfseries Case 1.} $i\neq s$, $j\neq t$.

Since $H$ is connected, there exists a path $i=i_1e_1i_2\cdots i_{\ell}e_{\ell}i_{\ell+1}=s$.
Since $r\geqslant 3$, there exist $j'$ such that $j'\neq j$, $j'\neq t$. From the definition
of $\widetilde{H}$, if $\ell$ is odd, we have
\[
\begin{cases}
(i_h,j)\sim (i_{h+1},j'), & h=1,3,\ldots,\ell-2,\\
(i_k,j')\sim (i_{k+1},j), & k=2,4,\ldots,\ell-1,\\
(i_{\ell},j)\sim (s,t).
\end{cases}
\]
If $\ell$ is even, we obtain
\[
\begin{cases}
(i_h,j)\sim (i_{h+1},j'), & h=1,3,\ldots,\ell-1,\\
(i_k,j')\sim (i_{k+1},j), & k=2,4,\ldots,\ell-2,\\
(i_{\ell},j')\sim (s,t).
\end{cases}
\]
Hence there exists a walk connecting $(i,j)$ and $(s,t)$.

\noindent{\bf Case 2.} $i=s$, $j\neq t$.

Since $r\geqslant 3$, there exist $i'$ and $j'$ such that $i'\neq i$, $j'\neq j$,
$j'\neq t$. According to Case 1 we know that there is a path connecting $(i,j)$ and
$(i',j')$. Noting that $i'\neq s$ and $j'\neq t$, there is a path connecting
$(i',j')$ and $(s,t)$ by Case 1. So there exists a walk connecting $(i,j)$ and $(s,t)$,
as desired. The proof of the claim is completed.

\noindent{\bf Claim 2.}\,  $\mathcal{A}(\widetilde{H})=\mathcal{A}(H)\otimes\mathcal{B}$.

\noindent{\bf Proof of Claim 2.} From the definition of $\widetilde{H}$, it follows that
\[
(\mathcal{A}(\widetilde{H}))_{(i_1,j_1)(i_2,j_2)\cdots(i_r,j_r)}=
\begin{cases}
\frac{1}{(r-1)!}, & \text{if}~\{i_1,i_2,\ldots,i_r\}\in E(H),b_{j_1j_2\cdots j_r}=1,\\
0, & \text{otherwise}.
\end{cases}
\]
According to Definition \ref{defn:Direct product}, $\mathcal{A}(H)\otimes\mathcal{B}$ is an order $r$
and dimension $rn$ tensor, whose entries are given by
\[
(\mathcal{A}(H)\otimes\mathcal{B})_{(i_1,j_1)(i_2,j_2)\cdots(i_r,j_r)}
=a_{i_1i_2\cdots i_r}b_{j_1j_2\cdots j_r}.
\]
If $\{i_1,i_2,\ldots,i_r\}\in E(H)$ and $b_{j_1j_2\cdots j_r}=1$, then
\[
(\mathcal{A}(H)\otimes\mathcal{B})_{(i_1,j_1)(i_2,j_2)\cdots(i_r,j_r)}=\frac{1}{(r-1)!},
\]
and $0$ otherwise. Hence
$\mathcal{A}(\widetilde{H})=\mathcal{A}(H)\otimes\mathcal{B}$, as desired.

\noindent{\bfseries Claim 3.}~~$\rho(\mathcal{B})=(r-1)!$.

\noindent{\bf Proof of Claim 3.}
Let $\bm{e}=(1,1,\ldots,1)^{\mathrm{T}}\in\mathbb{R}^r$.
It follows from Theorem \ref{relaigh} that
\[
\rho(\mathcal{B})\geqslant\frac{\bm{e}^{\mathrm{T}}(\mathcal{B}\bm{e})}{||\bm{e}||_r^r}
=\frac{r!}{r}=(r-1)!.
\]
On the other hand, let $\bm{z}=(z_1,z_2,\ldots,z_r)\in\mathbb{R}^r_{+}$ be a nonnegative
eigenvector corresponding to $\rho(\mathcal{B})$ with $||\bm{z}||_r=1$. By AM-GM inequality,
we have
\[
\rho(\mathcal{B})=\bm{z}^{\mathrm{T}}(\mathcal{B}\bm{z})=r!z_1z_2\cdots z_r\leqslant
r!\left(\frac{z_1^r+z_2^r+\cdots +z_r^r}{r}\right)=(r-1)!,
\]
with equality holds if and only if
\[
z_1=z_2=\cdots=z_r=\frac{1}{\sqrt[r]{r}}.
\]
Therefore,  $\rho(\mathcal{B})=(r-1)!$.
The proof of the claim is completed.

\noindent{\bf Claim 4.}~~$\rho(\widetilde{H})=(r-1)!\rho(H)$.

\noindent{\bf Proof of Claim 4.}
By Claim 2, $\mathcal{A}(\widetilde{H})=\mathcal{A}(H)\otimes\mathcal{B}$.
We consider the following two cases depending on whether or not $H$ is connected.

\noindent{\bf Case 1}. $H$ is connected.

Since $H$ is connected, we have $\mathcal{A}(H)$ is weakly irreducible.
From Theorem \ref{thm:Perron-Frobenius}, let $\bm{u}$ be the positive eigenvector
corresponding to the eigenvalue $\rho(H)$. Then by Theorem \ref{thm:Spectral H times G},
$\rho(H)\rho(\mathcal{B})$ is an eigenvalue of $\mathcal{A}(H)\otimes\mathcal{B}$
with a positive eigenvector $\bm{u}\otimes \bm{e}$. By Claim 1, $\widetilde{H}$
is connected. It follows from Theorem \ref{thm:Perron-Frobenius} and Claim 3 that $(r-1)!\rho(H)$
must be the spectral radius of $\mathcal{A}(H)\otimes\mathcal{B}$,
i.e., $\rho(\widetilde{H})=(r-1)!\rho(H)$.

\noindent{\bfseries Case 2}. $H$ is disconnected.

Let $\varepsilon>0$ and $\mathcal{A}_{\varepsilon}=\mathcal{A}(H)+\varepsilon\mathcal{J}_1$,
$\mathcal{B}_{\varepsilon}=\mathcal{B}+\varepsilon\mathcal{J}_2$, where $\mathcal{J}_1$ and
$\mathcal{J}_2$ are order $r$ tensors with all entries $1$ with dimension $n$ and $r$,
respectively. Then $\mathcal{A}_{\varepsilon}$ and $\mathcal{B}_{\varepsilon}$ are both
positive tensors, and therefore are weakly irreducible. Using the similar arguments as
Case 1, we have
\[
\rho(\mathcal{A}_{\varepsilon}\otimes\mathcal{B}_{\varepsilon})=
\rho(\mathcal{A}_{\varepsilon})\rho(\mathcal{B}_{\varepsilon}).
\]
Notice that the maximal absolute value of the roots of a complex polynomial is a continuous
function on the coefficients of the polynomial. Take the limit $\varepsilon\to 0$ on both
sides of the above equation, we obtain the desired result. The proof of the claim is completed.

It is clear that $\widetilde{H}$ is an $r$-partite hypergraph with partition
\[
V(\widetilde{H})=\bigcup_{i=1}^r\left(V(H)\times\{i\}\right).
\]
We define a vector $\bm{x}\in\mathbb{R}^{rn}$ as follows:
\begin{equation}
\label{eq:x}
x_{(i,j)}=\begin{dcases}
\frac{a_i}{\sqrt[r]{rn}}, & \text {if}~i=1,2,\ldots,n,j=1,\\
\frac{1}{\sqrt[r]{rn}}, & \text{otherwise},
\end{dcases}
\end{equation}
where $a_1$, $a_2$, $\ldots$, $a_n\geqslant 0$ and $a_1^r+a_2^r+\cdots+a_n^r=n$.
It is obvious that $d_{\widetilde{H}}((i,j))=(r-1)!d_i$ for any $i\in[n]$, $j\in[r]$.
By Theorem \ref{relaigh}, we deduce that
\begin{equation}
\begin{aligned}\label{eq:rho(tilde(H))}
\rho(\widetilde{H})\geqslant \bm{x}^{\mathrm{T}}(\mathcal{A}(\widetilde{H})\bm{x})
& =r\sum_{\{(i_1,j_1),\ldots,(i_r,j_r)\}\in E(\widetilde{H})}
x_{(i_1,j_1)}x_{(i_2,j_2)}\cdots x_{(i_r,j_r)}\\
& =r\left[\sum_{i=1}^n\frac{a_i}{\sqrt[r]{rn}}\cdot
\left(\frac{1}{\sqrt[r]{rn}}\right)^{r-1}\cdot (r-1)!d_i\right]\\
& =\frac{(r-1)!}{n}\sum_{i=1}^na_id_i.
\end{aligned}
\end{equation}
It follows from $a_1^r+a_2^r+\cdots+a_n^r=n$ and H\"older inequality that
\begin{align*}
\sum_{i=1}^na_id_i & \leqslant\left(\sum_{i=1}^na_i^r\right)^{\frac{1}{r}}\cdot
\left(\sum_{i=1}^nd_i^{\frac{r}{r-1}}\right)^{\frac{r-1}{r}}\\
& =\sqrt[r]{n}\left(\sum_{i=1}^nd_i^{\frac{r}{r-1}}\right)^{\frac{r-1}{r}}
\end{align*}
with equality  if and only if
\begin{equation}\label{eq:value of ai}
a_i=\frac{\sqrt[r]{n}d_i^{\frac{1}{r-1}}}{\sqrt[\leftroot{-2}\uproot{6}r]{\sum_{i=1}^nd_i^{\frac{r}{r-1}}}},~
i=1,2,\ldots,n.
\end{equation}
Now we set $a_i$ as \eqref{eq:value of ai}. In the light of \eqref{eq:rho(tilde(H))} and
\eqref{eq:value of ai} we have
\begin{equation}
\label{eq:Equality hold}
\rho(H)=\frac{\rho(\widetilde{H})}{(r-1)!}\geqslant
\frac{1}{n}\sum_{i=1}^na_id_i
=\left(\frac1n\sum_{i=1}^nd^{\frac{r}{r-1}}_i\right)^{\frac{r-1}{r}}.
\end{equation}

Now we give a characterization of extremal hypergraphs achieving the equality in \eqref{eq:Equality hold}.
Suppose first the equality holds in \eqref{eq:Equality hold}. Then the vector
$\bm{x}\in\mathbb{R}^{rn}$ defined by \eqref{eq:x} is an eigenvector corresponding to
$\rho(\widetilde{H})$ by Theorem \ref{relaigh}. Note that $H$ is connected, by
Theorem \ref{thm:Perron-Frobenius}, we let $\bm{u}=(u_1,u_2,\ldots,u_n)^{\mathrm{T}}\in\mathbb{R}_{++}^n$
be a positive eigenvector corresponding to $\rho(H)$. From Theorem \ref{thm:Spectral H times G}
and Claim 4, it follows that $\bm{u}\otimes\bm{e}$ is a positive eigenvector to $\rho(\widetilde{H})$.
By Claim 1 and Theorem \ref{thm:Perron-Frobenius}, $\widetilde{H}$ is connected, and we see that $\bm{x}$ and $\bm{u}\otimes\bm{e}$ are
linear dependence. Notice that
\[
\bm{x}=\frac{1}{\sqrt[r]{rn}}(\underbrace{a_1,1,\ldots,1}_r,\underbrace{a_2,1,\ldots,1}_r,\ldots,
\underbrace{a_n,1,\ldots,1}_r)^{\mathrm{T}}\in\mathbb{R}_{++}^{rn}
\]
and
\[
\bm{u}\otimes\bm{e}=(\underbrace{u_1,u_1,\ldots,u_1}_r,
\underbrace{u_2,u_2,\ldots,u_2}_r,\ldots,
\underbrace{u_n,u_n,\ldots,u_n}_r)^{\mathrm{T}}\in\mathbb{R}_{++}^{rn}.
\]
Consequently, $u_1=u_2=\cdots=u_n$, which implies that $H$ is regular.

Conversely, if $H$ is a  connected regular hypergraph, it is easy to see that the equality \eqref{eq:Equality hold} holds.
\qed

\begin{remark}
It is known that the spectral radius $\rho(H)$ of $H$ is greater than or equal to
the average degree \cite{Cooper:Spectra Uniform Hypergraphs}, i.e.,
\[
\rho(H)\geqslant\frac{\sum_{i=1}^nd_i}{n}.
\]
It follows from PM inequality that
\[
\left(\frac1n\sum_{i=1}^nd^{\frac{r}{r-1}}_i\right)^{\frac{r-1}{r}}
\geqslant\frac{\sum_{i=1}^nd_i}{n}.
\]
Therefore Theorem \ref{thm:Main result} has a better estimation for spectral radius of $H$.
\end{remark}

\section{Proof of Theorem 1.2}
In this section, we shall give a proof of Theorem \ref{Signless Laplacian}.

\noindent{\bfseries Proof of Theorem \ref{Signless Laplacian}.}\,
Let $H$ be a connected $r$-uniform hypergraph with vertex set $V(H)=[n]$.
Let $\widetilde{H}$ be the $r$-uniform hypergraph as defined in
Theorem \ref{thm:Main result}. Suppose that $\mathcal{B}$ is the order $r$ and dimension $r$
tensor given by \eqref{eq:B}, and $\mathcal{I}_r$ is the unit tensor of order $r$ and
dimension $r$. We have the following claims.

\noindent{\bfseries Claim 5.}
$\mathcal{Q}(\widetilde{H})=
(r-1)!(\mathcal{D}(H)\otimes\mathcal{I}_r)+\mathcal{A}(H)\otimes\mathcal{B}$.

\noindent{\bfseries Proof of Claim 5.}
Recall that $d_{\widetilde{H}}((i,j))=(r-1)!d_i$ for any $i\in[n]$, $j\in[r]$, we have
\[
\mathcal{D}(\widetilde{H})_{(i_1,j_1)(i_2,j_2)\cdots(i_r,j_r)}=
\begin{cases}
(r-1)!d_i, & \text{if}~i_1=\cdots=i_r=i,j_1=\cdots=j_r=j,\\
0, & \text{otherwise}.
\end{cases}
\]
On the other hand, by Definition \ref{defn:Direct product} we obtain that
\[
[(r-1)!(\mathcal{D}(H)\otimes\mathcal{I}_r)]_{(i_1,j_1)(i_2,j_2)\cdots(i_r,j_r)}=
(r-1)!\mathcal{D}(H)_{i_1i_2\cdots i_r}(\mathcal{I}_r)_{j_1j_2\cdots j_r}.
\]
If $i_1=i_2=\cdots=i_r=i$, $i\in[n]$ and $j_1=j_2=\cdots=j_r=j$, then
\[
[(r-1)!(\mathcal{D}(H)\otimes\mathcal{I}_r)]_{(i_1,j_1)(i_2,j_2)\cdots(i_r,j_r)}=(r-1)!d_i
\]
and $0$ otherwise. It follows that
\[
\mathcal{D}(\widetilde{H})=(r-1)!(\mathcal{D}(H)\otimes\mathcal{I}_r).
\]
Therefore, we have
\[
\mathcal{Q}(\widetilde{H})=\mathcal{D}(\widetilde{H})+\mathcal{A}(\widetilde{H})
=(r-1)!(\mathcal{D}(H)\otimes\mathcal{I}_r)+\mathcal{A}(H)\otimes\mathcal{B}.
\]
The proof of the claim is completed.

\noindent{\bfseries Claim 6.} $\rho(\mathcal{Q}(\widetilde{H}))=(r-1)!\rho(\mathcal{Q}(H))$.

\noindent{\bfseries Proof of Claim 6.}
Since $H$ is connected,  $\mathcal{Q}(H)$ is weakly irreducible. From
Theorem \ref{thm:Perron-Frobenius}, we let $\bm{u}$ be the positive eigenvector to
$\rho(\mathcal{Q}(H))$. Let $\bm{e}=(1,1,\ldots,1)^{\mathrm{T}}\in\mathbb{R}^r_{++}$.
By Proposition \ref{prop}, we deduce that
\begin{align*}
\mathcal{Q}(\widetilde{H})(\bm{u}\otimes\bm{e}) & =
[(r-1)!(\mathcal{D}(H)\otimes\mathcal{I}_r)+
\mathcal{A}(H)\otimes\mathcal{B}](\bm{u}\otimes\bm{e})\\
& =[(r-1)!(\mathcal{D}(H)\otimes\mathcal{I}_r)](\bm{u}\otimes\bm{e})+
[\mathcal{A}(H)\otimes\mathcal{B}](\bm{u}\otimes\bm{e})\\
& =(r-1)!(\mathcal{D}(H)\bm{u})\otimes(\mathcal{I}_r\bm{e})+
(\mathcal{A}(H)\bm{u})\otimes(\mathcal{B}\bm{e})\\
& =(r-1)!(\mathcal{D}(H)\bm{u})\otimes\bm{e}+
(r-1)![(\mathcal{A}(H)\bm{u})\otimes\bm{e}]\\
& =(r-1)![\mathcal{D}(H)\bm{u}+\mathcal{A}(H)\bm{u}]\otimes\bm{e}\\
& =(r-1)!(\mathcal{Q}(H)\bm{u})\otimes\bm{e}.
\end{align*}
It follows from $\mathcal{Q}(H)\bm{u}=\rho(\mathcal{Q}(H))\bm{u}^{[r-1]}$ that
\[
\mathcal{Q}(\widetilde{H})(\bm{u}\otimes\bm{e})=
(r-1)!\rho(\mathcal{Q}(H))(\bm{u}\otimes\bm{e})^{[r-1]},
\]
which yields that $\bm{u}\otimes\bm{e}$ is a positive eigenvector of
$\mathcal{Q}(\widetilde{H})$ corresponding to $(r-1)!\rho(\mathcal{Q}(H))$.
Note that $H$ is connected, then $\widetilde{H}$ is connected by Claim 1. Therefore
$\mathcal{Q}(\widetilde{H})$ is weakly irreducible. By Theorem \ref{thm:Perron-Frobenius},
$(r-1)!\rho(\mathcal{Q}(H))$ is the spectral radius of signless Laplacian tensor
$\mathcal{Q}(\widetilde{H})$, as claimed.

Let $\bm{x}\in\mathbb{R}^{rn}$ be the column vector defined by \eqref{eq:x}.
By Theorem \ref{relaigh}, we have
\begin{equation}
\begin{aligned}\label{eq:rho(Q)}
\rho(\mathcal{Q}(\widetilde{H})) &
\geqslant\bm{x}^{\mathrm{T}}(\mathcal{Q}(\widetilde{H})\bm{x})
=\bm{x}^{\mathrm{T}}(\mathcal{D}(\widetilde{H})\bm{x})
+\bm{x}^{\mathrm{T}}(\mathcal{A}(\widetilde{H})\bm{x})\\
& =\sum_{j=1}^r\sum_{i=1}^n\mathcal{D}(\widetilde{H})_{(i,j)(i,j)\cdots(i,j)}x_{(i,j)}^r
+\bm{x}^{\mathrm{T}}(\mathcal{A}(\widetilde{H})\bm{x})\\
& =\frac{(r-1)!}{rn}\sum_{i=1}^n(a_i^r+r-1)d_i
+\bm{x}^{\mathrm{T}}(\mathcal{A}(\widetilde{H})\bm{x}).
\end{aligned}
\end{equation}
Furthermore, by AM-GM inequality, we  have
\begin{equation}
\label{eq:a_i^r+r-1}
a_i^r+r-1=a_i^r+\underbrace{1+1+\cdots+1}_{r-1}\geqslant ra_i,~i\in[n].
\end{equation}
So it follows from \eqref{eq:rho(tilde(H))}, \eqref{eq:rho(Q)} and \eqref{eq:a_i^r+r-1} that
\begin{align*}
\rho(\mathcal{Q}(\widetilde{H})) & \geqslant
\frac{(r-1)!}{rn}\sum_{i=1}^n(a_i^r+r-1)d_i
+\frac{(r-1)!}{n}\sum_{i=1}^na_id_i\\
& \geqslant\frac{2(r-1)!}{n}\sum_{i=1}^na_id_i.
\end{align*}
By \eqref{eq:value of ai}, we obtain
\begin{equation}
\label{eq:euqality}
\rho(\mathcal{Q}(H))=\frac{\rho(\mathcal{Q}(\widetilde{H}))}{(r-1)!}
\geqslant\frac{2}{n}\sum_{i=1}^na_id_i=
2\left(\frac1n\sum_{i=1}^nd^{\frac{r}{r-1}}_i\right)^{\frac{r-1}{r}}.
\end{equation}

Suppose that the equality holds in \eqref{eq:euqality}. Then the vector $\bm{x}\in\mathbb{R}^{rn}$
defined by \eqref{eq:x} is an eigenvector corresponding to $\rho(\mathcal{Q}(\widetilde{H}))$
by Theorem \ref{relaigh} and $a_i=1$, $i\in[n]$ by \eqref{eq:a_i^r+r-1}. Recall that
$\mathcal{Q}(\widetilde{H})$ is weakly irreducible. By Claim 6, $\bm{u}\otimes\bm{e}$
is a positive eigenvector to $\rho(\mathcal{Q}(\widetilde{H}))$. We see that $\bm{x}$
and $\bm{u}\otimes\bm{e}$ are linear dependence by Theorem \ref{thm:Perron-Frobenius}.
Therefore,  $d_1=d_2=\cdots=d_n$, which implies that $H$ is regular.

Conversely, if $H$ is a regular connected hypergraph, it is  straightforward to verify that the equality \eqref{eq:euqality} holds.
\qed

\begin{remark}
Recently, Nikiforov introduces the concept of odd-colorable hypergraphs \cite{Nikiforov2016},
which is a generalization of bipartite graphs. Let $r\geqslant2$ and $r$ be even. An $r$-uniform
hypergraph $H$ with $V(H)=[n]$ is called odd-colorable if there exists a map $\varphi:[n]\to [r]$
such that for any edge $\{i_1,i_2,\ldots,i_r\}$ of $H$, we have
\[
\varphi(i_1)+\varphi(i_2)+\cdots+\varphi(i_r)\equiv\frac{r}{2}\, (\mbox{\rm mod}\, r).
\]
It was proved that if $H$ is a connected $r$-uniform hypergraph, then
$\rho(\mathcal{L}(H))=\rho(\mathcal{Q}(H))$ if and only if $r$ is even
and $H$ is odd-colorable \cite{Yuan:Odd-colorable}. Therefore we have
\[
\rho(\mathcal{L}(H))\geqslant2\left(\frac1n\sum_{i=1}^nd^{\frac{r}{r-1}}_i\right)^{\frac{r-1}{r}}
\]
for a connected odd-colorable hypergraph $H$, which generalizes the result in \cite{AimeiYu}.
\end{remark}


\begin{thebibliography}{10}
\bibitem{Bretto}
A. Bretto, Hypergraph Theory: An Introduction, Springer, 2013.

\bibitem{K.C.Chang.etc:Perron-Frobenius Theorem}
K. C. Chang, K. Pearson, T. Zhang, Perron-Frobenius theorem for nonegative tensors,
Commun. Math. Sci. 6 (2008) 507--520.

\bibitem{Cooper:Spectra Uniform Hypergraphs}
J. Cooper, A. Dutle, Spectra of uniform hypergraphs, Linear Algebra Appl. 436 (2012) 3268--3299.

\bibitem{Nikiforov2010}
M. Fiedler, V. Nikiforov, Spectral radius and Hamiltonicity of graphs,
Linear Algebra Appl. 432 (2010) 2170--2173.

\bibitem{Friedland}
S. Friedland, A. Gaubert, L. Han, Perron-Frobenius theorems for nonnegative multilinear forms
and extensions, Linear Algebra Appl. 438 (2013) 738--749.

\bibitem{Hofmeister}
M. Hofmeister, Spectral radius and degree sequence, Math. Nachr. 139 (1988) 37--44.

\bibitem{Lim} L. H. Lim, Singular values and eigenvalues of tensors: a variational approach,
in: Proceedings of the IEEE International Workshop on Computational Advances in Multi-Sensor
Adaptive Processing, vol. 1, CAMSAP'05, 2005, pp. 129--132.

\bibitem{Nikiforov2006}
V. Nikiforov, Eigenvalues and degree deviation in graphs,
Linear Algebra Appl. 414 (2006) 347--360.

\bibitem{Nikiforov2007}
V. Nikiforov, Bounds on graph eigenvalues II,
Linear Algebra Appl. 427 (2007) 183--189.

\bibitem{Nikiforov}
V. Nikiforov, Analytic methods for uniform hypergraphs,
Linear Algebra Appl. 457 (2014) 455-535.

\bibitem{Nikiforov2016}
V. Nikiforov, Hypergraphs and hypermatrices with symmetric
spectrum, arXiv: 1605.00709v2, 2016.

\bibitem{Qi2005}
L. Qi, Eigenvalues of a real supersymmetric tensor,
J. Symbolic Comput. 40 (2005) 1302--1324.

\bibitem{Qi2013}
L. Qi, Symmetric nonnegative tensors and copositive tensors,
Linear Algebra Appl. 439 (2013) 228--238.

\bibitem{Qi2014}
L. Qi, $H^+$-eigenvalue of Laplacian and signless Laplacian tensors, Commun. Math. Sci.
12 (2014) 1045-1064.

\bibitem{Shao:General product}
J. Shao, A general product of tensors with applications,
Linear Algebra Appl. 439 (2013) 2350--2366.

\bibitem{Shao:Connected odd-bipartite}
J. Shao, H. Shan, B. Wu, Some spectral properties and characterizations of
connected odd-bipartite uniform hypergraphs, Linear Multilinear Algebra
63 (12) (2015) 2359--2372.

\bibitem{Yuan:Odd-colorable}
X. Yuan, L. Qi, J. Shao, O. Chen,
Some properties and applications of odd-colorable $r$-hypergraphs,
arXiv:1606.05045v4, 2016.

\bibitem{Yang:Nonegative Weakly Irreducible Tensors}
Y. Yang, Q. Yang, On some properties of nonegative weakly irreducible tensors,
arXiv: 1111.0713v2, 2011.

\bibitem{AimeiYu}
A. Yu, M. Lu, F. Tian, On the spectral radius of graphs,
Linear Algebra Appl. 387 (2004) 41--49.


\end{thebibliography}
\end{document}